\documentclass
[a4paper
,11pt
]{amsart}

\usepackage{amsthm,amsmath,amssymb}
\usepackage{extarrows}

\usepackage{tikz-cd}

\numberwithin{equation}{section}

\newtheorem{theorem}[equation]{Theorem}
\theoremstyle{remark}
\newtheorem{remark}[equation]{Remark}
\theoremstyle{definition}

\newtheorem{example}[equation]{Example}

\renewcommand{\epsilon}{\varepsilon}
\renewcommand{\phi}{\varphi}

\newcommand{\Alg}{\ccat{Alg}}

\newcommand{\Bord}{\ccat{Bord}}

\newcommand{\cat}{\mathcal}
\newcommand{\Cat}{\ccat{Cat}}

\newcommand{\ccat}{\mathrm}

\newcommand{\Cocor}{\ccat{Cocorr}}
\newcommand{\Com}{\ccat{Com}}

\newcommand{\ddeloop}{\field{\deloop}}
\newcommand{\deloop}{B}

\DeclareMathOperator{\End}{End}

\newcommand{\equivwith}{\simeq}

\newcommand{\field}{\mathbb}

\newcommand{\Fin}{\ccat{Fin}}

\newcommand{\Gpd}{\ccat{Gpd}}

\newcommand{\indents}{ }
\newcommand{\Init}{\ccat{Init}}

\newcommand*{\kokoni}[1]{\noindent\makebox[0pt][l]{#1}}
\newcommand{\kore}{\textbf}

\newcommand{\longto}{\longrightarrow}

\newcommand{\Multicat}{\ccat{Multicat}}

\newcommand{\op}{\mathrm{op}}

\newcommand{\oriented}{\mathrm{or}}

\newcommand{\Set}{\ccat{Set}}

\newcommand{\tensor}{\otimes}
\newcommand{\Tensor}{\bigotimes}

\newcommand*{\tsketas}[2]
{#1\makebox[0pt][l]{$\vphantom{#1}#2$}}

\newcommand{\unity}{\boldsymbol{1}}
\newcommand{\wasreteori}{\Theta}

\begin{document} 
\title[Coherence and categorified structure]{%
  Higher
  coherence and a generalization of
  higher categorified
  algebraic structures.}

\author{Matsuoka, Takuo}
\email{motogeomtop@gmail.com}

\begin{abstract}
  This article reflects (to some degree) a talk I gave in
  Hiroshima, and another in Himeji, in November 2016.
  I am grateful to the organizers of the symposia for the
  invitations.
\end{abstract}

\maketitle

\tableofcontents

\section{Introduction}
\label{sec:introduction}
\subsubsection{}

Discovery or recognition of the right kind of \emph{algebraic
structure} is
often important in the development of mathematical subjects.
Starting perhaps with ``group'' in Galois theory, a list of examples
would easily grow long.
One might wish to systematically find and treat algebraic structures,
especially for the use of \emph{higher categorical} ideas, since more various structures are available with higher
categorical dimensionality.
Importance of higher category theory comes for example, from the
necessity of it for both analysis and construction of
\emph{topological field theories} (TFT), as has
become clear from the proof and generalization of Baez and Dolan's
\emph{cobordism hypothesis} \cite{bdolan} by Lurie and Hopkins
\cite{tft}.
Higher category theory also seems promising since
\emph{categorification} (in the
sense of Crane \cite{crane-frenke,crane}, about which the most
influential pioneer may have been Grothendieck) has been a useful
method for finding important new structures.

In this article, we shall give a overview of the relatively
non-technical aspects of our work \cite{theory}, in the most
fundamental part of which, a
concrete understanding of \emph{higher categorical coherence} leads
to systematic views on some (quite general) kinds of algebraic
structure.
(On the other hand, we shall be able to touch on only a few of more
topics which seem
next fundamental, and are also treated in our work, building on what
was mentioned.
We refer the reader to \cite{theory} for those topics on which we
could not touch in this article.)
A central role will be played by a certain process resulting from
an inductivity in the structure of higher coherence
(Section~\ref{sec:coherence} below), which turns out to produce (for
good consequences) more general kinds of structure than the process of
categorification (or further, the process of \emph{lax}
categorification) does.
The subject seems fundamental so the consequences of further research
in the direction of our work seem likely to provide useful
new viewpoints and methods not simply in algebra or geometry (e.g.,
through algebro-geometric methodology), but also in various subjects
where algebraic methods are useful, such as
mathematics, science, engineering and so on.

\subsubsection{}
In order to proceed, we need a few reminders and notations.

\subsubsection{}

The notion of coloured operad/\emph{multicategory} is a
generalization\indents
of the notion of symmetric monoidal category in the sense that a
symmetric monoidal category $\cat{A}$ has underlying multicategory,
which we shall denote by $\wasreteori\cat{A}$, from which\indents
$\cat{A}$ can be recovered.
An object of $\wasreteori\cat{A}$ is an object of $\cat{A}$, and, for
a finite set $S$, an $S$-ary multimap $x=(x_s)_{s\in S}\to y$ in
$\wasreteori\cat{A}$ is a map $\Tensor_sx_s\to y$ in $\cat{A}$, where
$\tensor$ denotes the monoidal multiplication operation on
$\cat{A}$.
The structure of $\cat{A}$ is recovered from the structure of
$\wasreteori\cat{A}$ since we have a universal $S$-ary multimap
$x\to\Tensor_sx_s$ in $\wasreteori\cat{A}$ for every family $x$ of
objects of $\cat{A}$ indexed by $S$.

An important role of a multicategory is the role of \emph{controling
  algebras} over it.
For a multicategory $\cat{U}$, a \kore{$\cat{U}$-algebra} in a
symmetric monoidal category $\cat{A}$\indents is a functor
$\cat{U}\to\wasreteori\cat{A}$ of multicategories.
Thus, a multicategory is analogous to an \emph{algebraic theory}
in the sense of Lawvere \cite{lawvere}, which controls algebraic
structures of a specific kind in a more restricted context.

Many kinds of algebraic structure (in a symmetric monoidal
category)\indents
are indeed controlled by an uncoloured operad:
\begin{description}
\item[Commutative] by ``$\Com$'', the terminal multicategory
\item[Associative] by ``$E_1$'', the ``associative'' operad
\item[Bare object (no structure)] by ``$\Init$'', the initial
  uncoloured operad
\end{description}
and so on.
With colours, more various structures can be controlled.

\subsubsection{}

For a uncoloured operad $\cat{U}$, there is a notion of
$\cat{U}$-monoidal category, which \emph{categorifies} the notion of
$\cat{U}$-algebra in the
sense that a $\cat{U}$-monoidal category is morally a
``$\cat{U}$-algebra'' in categories.
Namely, a $\cat{U}$-monoidal structure is an analogue on a category,
of the structure of a $\cat{U}$-algebra.
By letting $\cat{U}$ vary, we obtain symmetric, associative, braided
\cite{joyal-street}, and other notions of monoidal category.

For a \emph{coloured} multicategory $\cat{U}$, a similar notion of
$\cat{U}$-monoidal category is also not difficult to define.
This is an important example for us, of a categorified kind of
structure.

\subsubsection{}

There are in fact, different kinds of multicategory, among which,
symmetric (discussed above), planar and braided \cite{fiedoro} are
quite commonly worked with.
Before getting back to the subject, let us give a description of
where those kinds come from, which will motivate development of a
general framework.

Given a symmetric multicategory $\cat{U}$, we obtain the relevant
notion, that of \emph{$\cat{U}$-graded} multicategory (where
\emph{grading} is in a sense generalizing that e.g., of a ring), from
the notion of $\cat{U}$-algebra through a
certain general process, which (as well as the output of which) we
call \emph{theorization} inspired by Lawvere's notion.
Examples include:
\begin{align*}
\text{$E_1$-graded}&=\text{planar}\\
\text{$E_2$-graded}&=\text{braided}\\
\text{$\Init$-graded multicategory}&=\text{category}.
\end{align*}
We consider a symmetric multicategory as \emph{ungraded}.
The notion of graded multicategory (similarly to that of symmetric
multicategory) comes in a form enriched in a symmetric monoidal
category.
It turns out that a $\cat{U}$-graded multicategory enriched in sets,
is an equivalent datum to a symmetric multicategory lying over
$\cat{U}$:
\[
  \Multicat_\cat{U}(\Set)=\Multicat(\Set)_{/\cat{U}}.
\]
(We may also replace $\Set$ with the groupoid-enrichd category $\Gpd$
of groupoids.
Namely, we have
\[
  \Multicat_\cat{U}(\Gpd)=\Multicat(\Gpd)_{/\cat{U}},
\]
which is rather natural in the case where $\cat{U}$ is enriched in
$\Gpd$, e.g., for $\cat{U}=E_2$.)

The process of theorization is quite general and can also start from
various other kinds of structure than $\cat{U}$-algebra.
In general, the theorized structure turns out to be a more general
kind than a categorified structure (and a lax categorified
structure).
For example, the notion of $\cat{U}$-graded multicategory is a common
generalization of the notions of multicategory and of (possibly
lax/oplax) $\cat{U}$-monoidal category.

\begin{example}
Consider a category $\cat{C}$ as a multicategory having only
\emph{unary} multimaps.
Then, over $\cat{C}$, the notions we have are as follows.
\begin{description}
\item[Algebra] Functor on $\cat{C}$ (``left $\cat{C}$-module'').
\item[Categorification] Functor $\cat{C}\to\Cat$.
\item[Theorization] In the case enriched in $\Set$, it is category
  $\cat{X}$ equipped with a functor $\cat{X}\to\cat{C}$, among which,
  categorifications (and their lax morphisms) correspond to
  \emph{fibrations}
  \[\begin{tikzcd}
      \tsketas{\cat{X}}{^\op}\arrow[d]\\
      \tsketas{\cat{C}}{^\op}
    \end{tikzcd}\]
  (and not necessarily Cartesian functors over $\cat{C}^\op$).
\end{description}
\end{example}

\subsubsection{}
It seems to be a reasonable guess that the notion of $\cat{U}$-graded
multicategory had been known before our work even in the form enriched
in a symmetric monoidal category, and if it had in the enriched form,
then it would
probably have been a folklore that the notion more generally gets
enriched in a ($\cat{U}\tensor E_1$)-monoidal category, namely, a
$\cat{U}$-monoidal object in the $2$-category of associative monoidal
categories, or equivalently, associative monoidal object in the
$2$-category of $\cat{U}$-monoidal categories.

The idea of theorization leads to an explanation of this enrichment
through a natural \emph{generalization}.
Indeed, the notion of $\cat{U}$-graded multicategory can be enriched
``along'' a ``$2$-dimensional'' algebraic structure of a
certain kind.
The kind of structure is what we call \kore{$\cat{U}$-graded
  $2$-theory}, and is a natural \emph{theorization} of the notion of
$\cat{U}$-graded multicategory.
(We refer to a multicategory also as a \kore{$1$-theory}.)
This enrichment exists by a general nature of theorization, but turns
out to specialize to the mentioned probable folklore.

We would like to take a glance at examples in the case
$\cat{U}=\Init$.
Note that a $2$-category is a categorified $\Init$-graded $1$-theory,
which therfore, is an instance of an $\Init$-graded \emph{$2$-theory}
by what we have remarked above.
We shall denote a $2$-category $\cat{C}$ by $\wasreteori\cat{C}$ when
we consider the former as an $\Init$-graded $2$-theory.

For an associative monoidal category $\cat{A}$, let $\deloop\cat{A}$
denote its ``categorical deloop''.
Thus, $\deloop\cat{A}$ is a $2$-category with a base object $*$,
which is the only object of $\deloop\cat{A}$, equipped with an
equivalence $\End_{\deloop\cat{A}}(*)\equivwith\cat{A}$ of
associative monoidal categories.

Our first example is about a generalization of delooping for planar
multicategories.

\begin{example}\label{ex:deloop}
For a planar multicategory $\cat{M}$, there is a functorial
construction of an $\Init$-graded $2$-theory $\ddeloop\cat{M}$, for
which we have
\[
  \ddeloop\wasreteori\cat{A}\equivwith\wasreteori\deloop\cat{A}
\]
for an associative monoidal category $\cat{A}$.
A category enriched\indents
in $\cat{M}$ is in fact an equivalent datum to an $\Init$-graded
$1$-theory enriched ``along'' $\ddeloop\cat{M}$.
\end{example}

For the next example, for a nicely behaving associative monoidal
category $\cat{A}$, let $\Alg_1(\cat{A})$ denote its ``Morita''
$2$-category (due to
B\'enabou \cite{benabou}) of associative algebras and bimodules.

\begin{example}\label{ex:morita-2-theory}
For a planar multicategory $\cat{M}$, there is a functorial
construction of an $\Init$-graded $2$-theory $\cat{A}lg_1(\cat{M})$,
for which we have
\[
  \cat{A}lg_1(\wasreteori\cat{A})\equivwith\wasreteori\Alg_1(\cat{A})
\]
for $\cat{A}$ as above.
There is forgetful functor
\[
  \cat{A}lg_1(\cat{M})\longto\ddeloop\cat{M}
\]
of $\Init$-graded $2$-theories, generalizing the forgetful lax functor
\[
  \Alg_1(\cat{A})\longto\deloop\cat{A}.
\]
A category enriched\indents
in $\cat{M}$, or equivalently, an $\Init$-graded $1$-theory enriched
along $\ddeloop\cat{M}$ (Example~\ref{ex:deloop}), has a canonical
lift\indents to a category enriched\indents
along $\cat{A}lg_1(\cat{M})$.
\end{example}

\begin{remark}
In the situation of Example~\ref{ex:morita-2-theory}, if $\cat{M}$ is
of the form $\wasreteori\cat{A}$, then there is another lift of the
enriched category resulting from the obvious section
$\deloop\cat{A}\to\Alg_1(\cat{A})$ to the forgetful lax functor,
which sends the base object of $\deloop\cat{A}$ to
the unit algebra $\unity$ in $\cat{A}$.
However, a multicategory does not in general have a unit object,
namely, the nullary tensor product.
\end{remark}

\subsubsection{}
We shall discuss the role of theorization in algebra, and a method for
actually theorizing various kinds of algebraic structures iteratively.

\section{Higher theories}
\label{sec:theory}

\subsubsection{}

The notion of theorization results from the desire to control
algeraic structures.
For example, ask what may control graded multicategories.
A traditional answer to this says that they are controlled by a
multicategory, at least under some restriction.
For example, for a mutlicategory $\cat{U}$, Baez and Dolan \cite{opetope}
constructs a multicategory $\cat{U}^+$, which they call the
\kore{slice} of $\cat{U}$, whose algebras in the category
of sets are $\cat{U}$-graded multicategories with restricted colours:
\[
  \Alg_{\cat{U}^+}(\Set)
  =\Multicat_{\text{Same colours as $\cat{U}$'s}}(\Set)_{/\cat{U}}.
\]
For example, in the case $\cat{U}=E_1$, the slice multicategory
controls uncoloured planar operads.

While slicing was the main construction\indents
for Baez and Dolan, a much simpler structure can in fact control all
$\cat{U}$-graded multicateogries without a restriction on colours.
Namely, they are controlled\indents by a certain (ungraded)
\emph{$2$-theory} $\wasreteori\cat{U}$:
\[
  \Multicat_\cat{U}=\Alg_{\wasreteori\cat{U}},
\]
where the construction of $\wasreteori\cat{U}$ is in a simple one step
(of replacing each composition operation in $\cat{U}$\indents
with the ``bimodule'' ``corepresented'' by it \cite{theory}), which is
essentially as simple as the construction of the underlying
multicategory $\wasreteori\cat{A}$ (recalled in
Section~\ref{sec:introduction}) of a symmetric monoidal category
$\cat{A}$.

The key for this simplicity (as well as the ability to treat colours
naturally) is that a $2$-theory is a \emph{$2$-dimensional}
structure.
It appears
\begin{quote}
  \emph{natural to see $1$-dimensional structure
    (such as multicategory)\indents
  as algebra over a $2$-dimensional structure.}
\end{quote}

\subsubsection{}
There are further technical advantages in introducing higher
dimensional structures rather than staying in the world of
multicategories.
For example, the purpose of Baez and Dolan was to give a definition
of an $n$-category, which took a few more steps from the slice
construction.
On the other hand, we can further theorize the notion of $2$-theory
iteratively with the same method (to be described in
Section~\ref{sec:coherence}) and the resulting ``\kore{higher
  theories}'' or more specifically, ``\kore{$n$-theories}'', already
include $n$-categories, which are iterated categorifications.
Thus, (iterative) theorization is a more direct route to
$n$-categories\indents
than the slicing.

Dimensionality is the key here again: an $n$-category is
naturally \emph{$n$-dimensional} as is an $n$-theory, rather than
$1$-dimensional like a multicategory.

As we have seen in Section~\ref{sec:introduction}, theorization also
helps with \emph{enriching} notions.

\subsubsection{}
One may wonder what examples of higher theories there are.
Higher theories arise for example, through various general
constructions\indents
from various inputs.
\begin{itemize}
\item ``Delooping'' construction $\ddeloop$, of an ($n+1$)-theory
  from an $n$-theory.
  For a symmetric monoidal category $\cat{A}$, which one may call a
  categorified ``$0$-theory'', there is an equivalence
  \[
    \ddeloop^n\wasreteori\cat{A}\equivwith\wasreteori^{n+1}\deloop^n\cat{A}
  \]
  of ($n+1$)-theories, where the right hand side is the
  ($n+1$)-theory which corresponds to the symmetric monoidal
  ($n+1$)-category $\deloop^n\cat{A}$, the $n$-fold deloop of
  $\cat{A}$.
\item The ``pull-back'' of grading (Section~\ref{sec:grade}), and
  ``push-forward''
  constructions which respectively have a universal property as the
  left and the right adjoint.
\item A generalization of the Day convolution \cite{day}.
\end{itemize}
Many of these constructions \emph{raise} theoretic order, so
interests at least in multicategories, i.e., $1$-theories, would
naturally lead eventually to interests in all higher theories.

There are also some concrete constructions\indents
of higher theories in more specific situations.

\section{Meaning of theorization}

\subsubsection{}

The effectiveness of use of multicategories results from the ability
to express algebras of interest as functors of multicategories
(whenever possible).
This is analogous to a lesson learned in linear algebra: matrices and
their multiplication become most clearly understandable when one
regards matrices as linear maps between vector spaces.

Theorization is a method for finding for a given kind of algebraic
structure, say ``X-algebra'', a new kind of structure, say
``X-theory'', such that structures similar (in a way) to X-algebras
can be expressed as (``coloured'') morphisms of X-theories.
(We can naturally include colours when an X-algebra may have colours
as a multicategory does for example.)
We consider structures expressible as coloured morphisms of X-theories
as
controlled by the source X-theory, and enriched ``along'' the target
X-theory.

The relation of theorization with categorification is understood from
this point of view.
For us, relevance of categorification comes from the principle that
\begin{quote}
  \emph{a given notion ``X-algebra'', makes sense\indents
    in a categorified form\indents
    of X-algebra.}
\end{quote}
E.g., for a multicategory $\cat{U}$, there is a notion of
$\cat{U}$-algebra in a \emph{$\cat{U}$-monoidal} category
$\cat{A}$.\indents
Such a structure is equivalent to (or
defined as) a lax $\cat{U}$-monoidal functor
$\unity^0_\cat{U}\to\cat{A}$, where $\unity^0_\cat{U}$ denotes the
unit $\cat{U}$-monoidal category.
This is an \emph{enriched} notion\indents
of $\cat{U}$-algebra, and in general, a theorized structure
($\cat{U}$-graded multicategory in this case) is a more general place
where the original notion (the notion of $\cat{U}$-algebra here) can
be enriched.

\subsubsection{}

To be more specific, by theorizing a notion ``X-algebra'' to a notion
``X-theory'', we want a consequence of the form ``for a categorified
form $\cat{A}$ of $X$-algebra,
\begin{quote}
  \emph{the datum of an X-algebra in $\cat{A}$\indents
    is equivalent to the datum of a ``coloured'' functor
    $\unity^\mathrm{Theory}_\mathrm{X}\to\wasreteori\cat{A}$ of
    X-theories}'',
\end{quote}
where
\begin{itemize}
\item $\unity^\text{Theory}_\mathrm{X}$ denotes the terminal
  X-theory,
\item $\wasreteori\cat{A}$ denotes $\cat{A}$ \emph{as} an X-theory,
\item we omit explanation of colouring, which is only technical.
\end{itemize}
By an X-algebra in $\cat{A}$,\indents we mean a suitably coloured lax
functor $\unity^\text{Algebra}_X\to\cat{A}$ of $X$-algebras.

Once we have a notion of X-theory satisfying this, it is reasonable to
define \kore{algebra} over an X-theory $\cat{U}$ as
a (coloured) functor on $\cat{U}$, so we will have the picture
\[\begin{tikzcd}
\text{X-theory}
&\ni
&\cat{U}\arrow{d}{\text{control}}
&\tsketas{\unity}{^\text{Theory}_\mathrm{X}}
\arrow{d}{\text{control}}\\
\text{X-algebra}\arrow[maps
to,dashed]{rr}{\substack{\text{``grade'' by $\cat{U}$}\\
    \text{(see below)}}}
\arrow[maps to,dashed]{u}{\text{theorize}}
&&\text{$\cat{U}$-algebra}&\text{X-algebra\kokoni{.}}
\end{tikzcd}\]

We will obtain enriched notions of algebra, which are interrelated as
follows.
\[\begin{tikzcd}[column sep=large]
\substack{\displaystyle\text{Functor}\\
  \displaystyle\cat{U}\to\cat{V}}
&\substack{\displaystyle\text{$\cat{U}$-algebra}\\
  \displaystyle\text{in $\cat{V}$}}
\arrow[maps to]{l}[swap]{\substack{\text{forbid}\\
    \text{colours}}}
\arrow[maps to]{d}[swap]{\cat{V}=\wasreteori\cat{A}}
\arrow[maps to]{r}{\cat{U}=\unity^\text{Theory}_\mathrm{X}}
&\substack{\displaystyle\text{X-algebra}\\
  \displaystyle\text{in $\cat{V}$}}
\arrow[maps to]{d}[swap]{\cat{V}=\wasreteori\cat{A}}\\
&\substack{\displaystyle\text{$\cat{U}$-algebra}\\
  \displaystyle\text{in $\cat{A}$}}
\arrow[maps to]{r}[swap]{\cat{U}=\unity^\text{Theory}_\mathrm{X}}
&\substack{\displaystyle\text{X-algebra}\\
  \displaystyle\text{in $\cat{A}$}}
\arrow[maps to]{r}{\substack{\text{forbid}\\
    \text{colours}}}
&\substack{\displaystyle\text{Lax functor}\\
  \displaystyle\unity^\text{Algebra}_\mathrm{X}\to\cat{A},}
\end{tikzcd}\]
where
\begin{description}
\item[$\cat{U}$, $\cat{V}$] X-theories
\item[$\cat{A}$] categorified X-algebra.
\end{description}

\subsubsection{}
The way how we can actually theorize in the desired manner in many
situations, will be discussed in Section~\ref{sec:coherence}.

\section{Graded higher theories}
\label{sec:grade}

\subsubsection{}

Let us describe a large system to which iterative theorization of the
notion of multicategory leads.

\subsubsection{}
Given an $n$-theory $\cat{U}$, there is a notion of
\kore{$\cat{U}$-graded $n$-theory}, which theorizes the notion of
$\cat{U}$-algebra.

The following theorem generalizes our method of controlling graded
multicategories mentioned in Section~\ref{sec:theory}.
As before, it is easy to construct an ($n+1$)-theory
$\wasreteori\cat{U}$ ``underlying'' $\cat{U}$.

\begin{theorem}
The following forms of datum are equivalent:
\begin{itemize}
\item a $\cat{U}$-graded $n$-theory
\item a $\wasreteori\cat{U}$-algebra.
\end{itemize}
\end{theorem}

There is also a notion of \kore{$\cat{U}$-graded ($n+1$)-theory},
which naturally theorizes the notion of $\cat{U}$-graded $n$-theory,
and we have the following theorem.
\begin{theorem}
The following forms of datum are equivalent:
\begin{itemize}
\item a $\cat{U}$-graded ($n+1$)-theory
\item a $\wasreteori\cat{U}$-graded ($n+1$)-theory.
\end{itemize}
\end{theorem}

By the general principle, a $\cat{U}$-graded ($n+1$)-theory is a
natural place where the notion of $\cat{U}$-graded $n$-theory can be
enriched.
In the case $n=1$, this specializes to the enrichment\indents
of the notion of $\cat{U}$-graded multicategory\indents
in a ($\cat{U}\tensor E_1$)-monoidal category, as mentioned in
Section~\ref{sec:introduction}.

The above theorems lead to iterative theorization.
Namely, they make it immediate to figure out the right notion of
\emph{$\cat{U}$-graded $m$-theory} for $m\ge n+2$.

\subsubsection{}

There are also graded ``lower'' theories.

Let $\cat{U}$ be an $n$-theory.
Then, for an integer $m$ such that $0\le m\le n-1$, there is a notion
of \emph{$\cat{U}$-graded $m$-theory} so that
\begin{quote}
  the notion of $\cat{U}$-graded ($m+1$)-theory theorizes the notion
  of $\cat{U}$-graded $m$-theory,
\end{quote}
and
\begin{quote}
  the notion of $\cat{U}$-algebra is equivalent to the notion of
  $\cat{U}$-graded ($n-1$)-theory.
\end{quote}

Moreover, there is a notion of $\ell$-theory graded by a
\emph{$\cat{U}$-graded} $m$-theory.
We obtain for example, equivalence of the following notions for every
integer $m\ge 0$:
\begin{itemize}
\item $\unity^m_\cat{U}$-graded $\ell$-theory
\item $\cat{U}$-graded $\ell$-theory
\end{itemize}
where $\unity^m_\cat{U}$ denotes the terminal $\cat{U}$-graded
$m$-theory.

\section{More general higher theories: some examples}
\label{sec:general-theory}
\subsubsection{}
We can also theorize some structures\indents
which (unlike a multicategory) involve operations with multiple
inputs \emph{and} multiple outputs.
More precisely, consider each symmetric monoidal category $\cat{B}$
as \emph{controlling} symmetric monoidal functors
$\cat{B}\to\cat{A}$, where $\cat{A}$ varies among symmetric monoidal
categories.
This ``kind'' of structure ``controlled'' by $\cat{B}$, can be
theorized iteratively if $\cat{B}$ is generated in a certain nice
manner.
We refer to the resulting notion
as the notion of \kore{$\cat{B}$-graded} $n$-theory.

\begin{remark}
In this terminology, a symmetric multicategory will be a
``$\Fin$-graded'' $1$-theory, where $\Fin$ denotes the category of
finite sets made symmetric monoidal under the operations of disjoint
union.
Note the difference of this notion from a $\wasreteori\Fin$-graded
$1$-theory in the sense of Section~\ref{sec:introduction}.
In the terminology here, the term $1$-theory does not translate as
multicategory, which is the point of the generalization.
\end{remark}

For example, the notion of \emph{coloured properad} of Vallette
\cite{vallette} turns out\indents
to be equivalent to the notion of $1$-theory graded by $\Cocor(\Fin)$,
the cocorrespondence category (enriched in groupoids) on $\Fin$.
The notions of $\Cocor(\Fin)$-graded $n$-theory iteratively theorize
the notion of coloured properad.

\subsubsection{}
In another instance, we find structures of somewhat unexpected nature.
Namely, for an infinity $1$-category $\cat{C}$, there is a
``$\Bord^\oriented_1$-graded'' $1$-theory $\cat{Z}_\cat{C}$, where
$\Bord^\oriented_1$ denotes the oriented $1$-dimensional bordism
category (enriched in groupoids) ``controlling'' $1$-dimensional
oriented (equivalently, framed) TFT's, such that\indents
every ``$1$-dimensional TFT'' \emph{in $\cat{Z}_\cat{C}$}, namely,
functor
\[
  \unity^1_{\Bord_1}\longto\cat{Z}_\cat{C}
\]
of $\Bord^\oriented_1$-graded $1$-theories, is of the form
\[
  \unity^1_{\Bord_1}=\cat{Z}_{\unity}\xlongrightarrow{\cat{Z}_x}\cat{Z}_\cat{C}
\]
for an unique object
\[
  x\colon\unity\longto\cat{C}
\]
of $\cat{C}$.

This looks very different from field theories in the original,
\emph{untheorized} context, which, according to the cobordism
hypothesis, are classified by
\emph{dualizable} objects of a \emph{symmetric monoidal} infinity
$1$-category \cite{tft}.
Yet the notions are parallel, so a unified treatment (as far as that
goes) is at least natural.

\section{Higher coherence and iterative theorization}
\label{sec:coherence}
\subsubsection{}

The way how we can iteratively theorize\indents
kinds of structure discussed in the previous sections, is by using an
inductivity\indents
embedded in the structure of the
\emph{coherence} for higher associativity.

To describe the idea of this inductivity, suppose that we have the
following data.
\begin{description}
\item[$m$] a collection of operations wanting to be
  \emph{associative}, operating as maps in an symmeric monoidal
  infinity $1$-category $\cat{A}$.
  E.g., ``$m$ with a coherent associativity'' may\indents
  define an ``X-algebra'' enriched in $\cat{A}$.
\item[$m'$] collection of $2$-isomorphisms/homotopies giving an
  associativity of $m$.
\end{description}
Then we would wish to make $m'$ \emph{coherent}.

In the case where ``X-algebra'' (i.e., the kind of the structure
being constructed in the described situation) means any of the kinds
considered in the previous sections for theorization, there is a way
to organize $m'$ so we can see it as a collection of operations
themselves in
such a manner that \emph{associativity} of those \emph{operations}
$m'$ makes sense.
(This is explicitly visible for instance, in the definition of a
symmetric $n$-theory in \cite{theory}.)

\begin{remark}
It may not seem natural to try to see $m'$ as operations since we
required $m'$ to consist of $2$-isomorphisms rather than possibly
non-invertible $2$-morphisms.
However, for $m'$ in the suitably organized form as mentioned, we do
not actually need to require invertibility since, without
invertibility, we would still have the structure of an
(op)\emph{lax} X-algebra, with which we could be contented.
\end{remark}

Moreover, we obtain equivalence of the following forms of datum:
\begin{itemize}
\item a \emph{coherence} of the associativity $m'$ for $m$\indents
\item a coherent \emph{associativity} of $m'$ \emph{as operations}.
\end{itemize}

\subsubsection{}

The described inductivity is relevant to theorization since many
kinds of structure can be theorized using the idea that, if
$m$\indents
gives an ``X-algebra'' structure, then its theorization, ``X-theory'',
will be a kind\indents
defined (in the form enriched in a symmetric monoidal category
$\cat{A}$) by \emph{operations} $m'$\indents considered in the
``categorical deloop'' $\deloop\cat{A}$ (so
$2$-morphisms $m'$ are \emph{maps} in
$\cat{A}=\End_{\deloop\cat{A}}(*)$).

Since the structure is inductive, we can iterate theorization, and get
to a notion of $n$-theory for every integer $n\ge 0$,
generalizing in particular, the notion of \emph{$n$-category}.
The definition of an $n$-theory can be written explicitly
\cite{theory}.

\end{document}